\documentstyle[epsf,sectioneqno]{article}
\newtheorem{theorem}{Theorem}[section]
\newtheorem{lemma}[theorem]{Lemma}

\newtheorem{Algorithm}{Algorithm}[section]
  \title{    Approximate diagonalization in differential systems and an effective algorithm for     the computation of the spectral matrix  }
\author{    B.M.Brown, M.S.P. Eastham,D.K.R.M$^c$Cormack \\ 
 Department of Computer Science, Cardiff University of Wales, Cardiff, CF2 3XF, U.K. 
\and W.D.Evans \\ School of Mathematics, \\ University of Wales College of Cardiff, \\Senghennydd Road, Cardiff, CF2 4AG  U.K.  }
\date{}
\begin{document}
\maketitle


 \newcommand{\beq}{\begin{equation}}
\newcommand{\enq}{\end{equation}}
\newcommand{\s}{ {\cal{x}}}

\section{Introduction}
 In a recent paper  \cite{BEEKa94}, an extended Liouville-Green formula
\beq
y(x) =  \gamma_M(x) \{ 1 + \epsilon_M(x) \}\exp(\pm \int_X^xQ_M(t)dt)
 \label{eq:1.1}
\enq
was developed for solutions of the second-order differential 
equation
\beq
y''(x) - Q(x)y(x) = 0  \;\;\; (0 \leq x < \infty). \label{eq:1.2}
\enq
Here $ \gamma_M(x)  \sim Q^{-\frac{1}{4}}(x), \;Q_M(x)  \sim  Q^{\frac{1}{2}}(x)$ and $\epsilon_M(x) \rightarrow  0$ as $x \rightarrow  \infty,$ 
while $M
 (\geq 2)$ is an integer and $\gamma_ M $ and $Q_M$ can be defined in 
terms of $Q $ and its derivatives up to order $M - 1.$ The general 
form of   (\ref {eq:1.1})  had been obtained previously by  
Cassell \cite{JSC85} \cite{JSC86} \cite{JSC88} and Eastham \cite{MSPE87} \cite[ section 2.4] {MSPE89}.  In particular, the 
proof of   (\ref {eq:1.1})  in  \cite{MSPE87} and \cite{MSPE89} depended on the formulation of  (\ref{eq:1.2})   as 
a first-order system and then on a process of $M$ repeated 
diagonalizations of the coefficient matrices in a sequence of 
related differential
systems.
\par
   The main contribution of \cite{BEEKa94} to   (\ref {eq:1.1})  was to show that, for 
a general class of coefficients $Q,$ the magnitude of the error 
term $\epsilon_M(x) $ for large $x$ decreases as $M $ increases. This feature 
of $\epsilon_M(x) $ was then exploited in \cite{BEEKa94} for the case where $Q(x) = 
q(x) - \lambda$  and  (\ref{eq:1.2})   becomes the usual Sturm-Liouville equation 
with spectral parameter $\lambda$ . The smallness of $\epsilon_M(x) $ for a 
suitable choice of $M$ (such as $M = 6$) leads to an efficient 
numerical algorithm for estimating the Titchmarsh-Weyl function 
$m(\lambda )$ with precise global error bounds. In \cite{BEEKa94}, the case where 
$q(x) = -x^\alpha \; \; (0 < \alpha \leq 2)$ and 
  $ {\rm{Re}} \lambda \geq -1 $ was considered in detail.
\par
   In this paper, we develop these ideas for higher-order 
differential equations. Thus we consider the generalisation of 
  ( \ref {eq:1.1})  for the {\it{n}}-th order equation
\beq
                          y^{(n)}(x) - Q(x)y(x) = 0            \label{eq:1.3}   
\enq
again with improving estimates on $\epsilon_M$ as $M$ increases. In the 
case where $n$ is even, $n = 2\nu$, and $Q(x) = (-1)^\nu\{\lambda    - q(x)\} $, we 
have a   $\nu \times \nu$  spectral matrix $(m_{ij}(\lambda ))$ which corresponds to the 
Titchmarsh-Weyl function $m(\lambda )$ \cite{Eve64}, and we can again apply our 
estimates for $\epsilon_M$ to obtain numerical estimates for the spectral 
matrix.
 \par
There are however certain difficulties which arise when $ n > 
2$ and which were not present for the second-order equation 
 (\ref{eq:1.2})  . The first is that, when $n > 2,$ it is no longer possible 
to diagonalize the $n \times n$ matrices in the
associated first-order systems in the same explicit way as when 
$n = 2.$ We overcome this difficulty by adopting an approximate 
diagonalization process and this approach is discussed in sections 2 and 3. 
Second, this approximate procedure necessarily involves 
differential systems whose coefficient matrices contain a 
greater number of terms than when $n = 2$, and we have to develop 
a new   algorithm for generating and collating 
these terms. These matters are introduced in sections 4-5. Then in section 6, we   develop a computational algorithm
 to perform the repeated diagonalization. The 
estimation of error terms leading to $\epsilon_M$ presents fewer 
difficulties and it is also dealt with in section  6. The application 
to the spectral matrix  $(m_{ij}(\lambda ))$  is naturally less simple 
than when 
$n = 2$ because the underlying spectral theory involves several 
$L^2(0,\infty)$ solutions of the differential equation rather than just 
one such solution. We cover this application in sections 7 and 8 with special reference to the example $q(x)=-x^\alpha\; ( \alpha >0)$.
In section 9 we are able to give independent confirmation of 
our results when $\alpha =1$ in terms of the higher-order Airy equation. 
Finally, in section 10 we indicate possible extensions of our work and, 
in both sections 8 and 10, we comment on the effectiveness of our methods 
as compared with the recent alternative approach 
of Bennewitz et al. \cite{BBEMM95}.
\section{Diagonalization in differential systems}
   In this section we introduce the theoretical basis for 
estimating and improving error terms in the solution of 
differential systems. These error terms give rise to $\epsilon_M$ in 
 (\ref{eq:2.1})   and in the corresponding formulae for  (\ref{eq:1.3})   (see  (\ref{eq:3.14})   
in the next section). The main components of the discussion are 
 covered in the three subsections (a)-(c) below, and we frame 
the discussion in terms of the system
\beq
 Z'(x) = x ^\beta \{D + O(x^{-\gamma })\} Z(x ), \label{eq:2.1}             
\enq
where $D$ is a constant diagonal $n \times  n $ matrix,
\begin{displaymath}
D =  {\rm{dg}} (d_1, ..., d_n)
\end{displaymath}
with distinct $d_k$ and the $O-$term refers to $x  \rightarrow \infty $. 
Also, $ \beta > 0 $
and  $\gamma > 0$. If
in addition $ \beta -\gamma < -1 $,   (\ref{eq:2.1})   
has the standard Levinson form  \cite{NL48} \cite[section 1.3
 ]{MSPE89}
 \beq
 Z' = ( \Lambda + R)Z       \label{eq:2.2}              
\enq
in which  $\Lambda $ is diagonal and $R$ is $L(a, \infty )$, and the Levinson 
asymptotic theorem  \cite[  Theorem 1.3.1]{MSPE89} states that there are 
solutions $Z_k$ $(1  \leq k \leq  n)$ of  (\ref{eq:2.1})   such that
\beq
 Z_k(x) = (e_k +  \eta_k)\exp( \int_{a}^{x}t^\beta d_k dt) \label{eq:2.3}       
\enq
where $e_k$ is the unit coordinate vector in the $k-$direction and 
 $\eta_k   \rightarrow 0 $ as 
$x  \rightarrow \infty  $. The size of $ \eta_k$ is related to the size of 
$x^{\beta-\gamma } $ as $x  \rightarrow \infty   $
in a manner to be made precise in section 3.
\vspace{0.25in}
\par
\newcounter{rem1}
\par
(a) {\bf{\it{ Exact Diagonalization}}}  We first consider the effect of 
expressing the coefficient matrix in  (\ref{eq:2.1})   in its diagonal 
form and, as a start, we suppose that the O-term is simpley $x^{-\gamma}C$, where $C$ is a constant. Thus we write
\beq
                           T^{-1}(D + x^{-\gamma} C)T = D_1,       \label{eq:2.4}          
\enq
where $D_1$ is the diagonal matrix formed by the eigenvalues of $ D 
+ x^{- \gamma}C $ and $T$ is the approximate identity matrix formed by the 
eigenvectors. Then
\begin{displaymath}
       D_1 = D + O(x^{-\gamma} ), \;\;\; T = I + O(x^{-\gamma} ),\;\;\;  
       T' = O(x^{-\gamma -1})
\end{displaymath}
as $x   \rightarrow \infty$. The transformation
\beq
                                    Z = TW          \label{eq:2.5}            
\enq 
therefore takes  (\ref{eq:2.1})   into
\beq            
W' = x^\beta (D_1 - x^{-\beta}T^{-1}T') W = 
x^\beta \{D_1 + O(x^{-\beta-\gamma 
-1 })\} W.           \label{eq:2.6}
\enq 
Hence the effect of the transformation  (\ref{eq:2.5})   is to replace the 
perturbation
$x^{-\gamma}$  in  (\ref{eq:2.1})   by one of smaller magnitude $ x^{-\beta -\gamma -1} $
in  (\ref{eq:2.6}).
\par
   Even though  (\ref{eq:2.6})   does not have quite as simple a form as we 
specified in  (\ref{eq:2.1}), the argument which leads from  (\ref{eq:2.1})   to 
 (\ref{eq:2.6})   can again be applied to  (\ref{eq:2.6})   to yield a sequence of 
systems
\beq 
W^\prime_m = x^{\beta}(D_m + R_m)W_m  
\label{eq:2.7}                
\enq
with approximately constant diagonal matrices $D_m$ and 
perturbations $R_m$ of magnitude $x^{-m(\beta+1)- \gamma} \;\;\; 
(m = 1, 2, ...)$. The 
Levinson formula (2.3) can be applied to any system in the 
sequence, and the term  $\eta_k$ will have a corresponding order of 
magnitude. The transformation back to the original system  (\ref{eq:2.1})   
via the sequence of equations such as  (\ref{eq:2.5})   thus gives  successively more 
accurate asymptotic representations of the solutions of  (\ref{eq:2.1}). 
\par
The details of this process were given in \cite[ sections 2-3]{BEEKa94}  for the
case where  (\ref{eq:2.1})   arises from (1.2) and all the matrices are 
$2 \times 2$. These details depended on the explicit knowledge of all 
necessary eigenvalues and eigenvectors in equations such as 
 (\ref{eq:2.4})  . When $n > 2$, however, such explicit knowledge is not 
available in  (\ref{eq:2.4}). In this paper, we overcome this difficulty 
by modifying  (\ref{eq:2.4})   and considering instead an approximate 
diagonalization process which is based on the one introduced by 
Eastham \cite{MSPE86}
\cite[section 1.7 ]{MSPE89}.
\vspace{0.25in}
\par
(b) {\bf{\it{Approximate Diagonalization}}}  The matrix $T$ in  (\ref{eq:2.5})   has 
the form 
$T = I + P,$ where $P = O(x^{-\gamma} )$, and the idea now is to modify 
 (\ref{eq:2.5})   to
\beq
Z = (I + P)W    \label{eq:2.8}                 
\enq
with a different, but explicit $P$. It follows from \cite[ (1.7.2) 
and (1.6.13)]{MSPE89} that  (\ref{eq:2.8})   represents an approximation to  (\ref{eq:2.5})   
if $P$ is defined by
\beq                    
PD - DP = x^{-\gamma} (C -  {\rm{dg}}  C)  \label{eq:2.9}              
\enq
with $ {\rm{dg}}  P = 0.$ Thus the entries $p_{ij}$ in $P$ are defined by
\beq
p_{ij} = x^{-\gamma} c_{ij}/(d_j - d_i)   \;\;\;(i \neq   j).        
\label{eq:2.10}  
\enq
By  (\ref{eq:2.8})   and (2.9), we obtain in place of  (\ref{eq:2.6})  
\begin{eqnarray}
      W' &=& x^{\beta} \{ D + x^{-\gamma} (I + P)^{-1}(CP +  {\rm{dg}}  C) 
      - x^{-\beta}(I + P)^{-1}P'\} W \nonumber \\
         &=& x^{\beta} \{  D_1 + x^{-\gamma}  (I + P)^{-1}(CP - P {\rm{dg}}  C) 
         - x^{-\beta}(I + P)^{-1}P'\} W.\nonumber\\   \label{eq:2.11}
\end{eqnarray}  
Here 
\beq                        D_1 = D + x^{-\gamma}  {\rm{dg}}  C      \label{eq:2.12}       \enq
and, by  (\ref{eq:2.10})  , the other two groups of terms are respectively 
$O(x^{-2\gamma} )$ and  $O(x^{-\beta-\gamma-1})$, 
of which the latter dominates if  $\gamma > \beta 
+ 1.$ We note that  (\ref{eq:2.11})   contains more terms than  (\ref{eq:2.6}), but 
this is the price to be paid for having explicit terms derived 
from  (\ref{eq:2.10}).
\vspace{0.25in} 
\par
(c) {\bf{\it{Further approximation}}}  Repetition of the process  (\ref{eq:2.8})  -
 (\ref{eq:2.11}), but
with  (\ref{eq:2.11})   as the starting point, involves the diagonal 
entries of $D_1$ in place of $d_j$ and $d_i$ in  (\ref{eq:2.10}). Thus powers of $ x $
appear increasingly in the denominators of the $P-$matrices. This 
causes problems when we seek to develop computing and numerical 
procedures for implementing the repeated transformations, 
particularly because repeated differentiation of such matrices 
is involved, as indicated by the appearance of $P'$ in  (\ref{eq:2.11}). 
Accordingly, we shall retain the original constant entries $d_j$ 
and $d_i$ in the denominators as we go through the process. This 
further approximation adds new terms to the coefficient 
matrices  $\{...\}$  in the sequence of systems corresponding to 
 (\ref{eq:2.11}). It is this last method that we adopt in this paper and, 
after these introductory remarks, we defer further details to 
section 4 where we deal with the actual system  (\ref{eq:2.1})   which arises from 
 (\ref{eq:1.3}).
\section{Asymptotic formulae for solutions}
   The dominant form of the solutions of (1.3) for large $x$ is 
known subject to suitable conditions on $Q$ \cite[section 2.9 ]{MSPE89}. Our 
aim here is to obtain an explicit estimate for the error term 
in the asymptotic formulae, which is suitable for our purposes 
in  sections  4-6. We begin with the system formulation
\beq
                 Z' = ( \Lambda + R)Z  \label{eq:3.1}                    
\enq      
of (1.3) which is given in \cite[ section 2.9 ]{MSPE89}. Here
 \beq
              \Lambda             = Q^{1/n} {\rm{dg}} ( \omega_1, ...,\omega_n),               
  \label{eq:3.2}
\enq
where                       
\begin{displaymath}
\omega_k = \exp \{2(k - 1)\pi i/n \}
\end{displaymath}
are the $n-$th roots of unity, and
 \beq
                               R = -\Omega^{-1} \Lambda_1 \Omega,                     
\label{eq:3.3}
\enq
where                 
\beq\Lambda_1 = n^{-1}(Q'/Q) {\rm{dg}} (0, 1, ..., n - 1),        
                                               \label{eq:3.4}
\enq
while $\Omega$ and $\Omega^{-1}$ have $(j,k)$ entries $\omega^{ j-1}_k $
and $n^{-1}\omega_j^{-(k-1)} $
respectively. The connection between $Z$ and $y$ in  (\ref{eq:3.1})   and (\ref{eq:1.3}) 
is
\beq                   
Y =  \{ {\rm{dg}} (1, Q^{1/n}, Q^{2/n}, ..., Q^{(n-1)/n} )\}\Omega Z
\label{eq:3.5}
\enq
where $Y$ has components $y, y', ..., y^{(n-1)}$. It follows from 
 (\ref{eq:3.2})  - (\ref{eq:3.4})   that  (\ref{eq:3.1})   can be written in the form
\beq
                         Z' = Q^{1/n}(D + Q'Q^{-1-1/n}C)Z,            
\label{eq:3.6}                                                          
\enq
where $C$ is constant and
\beq
D =  {\rm{dg}} ( \omega_1, ..., \omega_ n).                
\label{eq:3.7}
\enq
In the case where $Q$ and its derivatives sufficiently resemble 
$x^\alpha \;\;(\alpha >0)$ and its derivatives, the system (\ref{eq:3.6})  essentially 
has the form (\ref{eq:2.1}) with 
$\gamma  = 1 + \alpha/n$, and the ideas in  section 2 for improving the error term 
apply. At this point, we note that the diagonal entries in $C$ 
are all $-(n - 1)(2n)^{-1}$
\cite[   (2.9.25)  ]{MSPE89} and this, together with (\ref{eq:2.12}), governs our 
definition of $D_M$ and $D_m$ in  (\ref{eq:3.10})   and (\ref{eq:4.3}) below.
\par
   The sequence of transformations to be defined in  section 4 takes 
(3.6) into the form
\beq
                           Z'_M = Q^{1/n}(D_M + R_M)Z_M,                
 \label{eq:3.8}    
\enq
where               
\beq
Z =  \{\prod_1^{M-1}  (I + P_m)\}Z_M=(I+P)Z_M,       \label{eq:3.9}         
\enq 
\beq
             D_M = D - (n - 1)(2n)^{-1}Q'Q^{-1-1/n}I +  \Delta_M       
             \label{eq:3.10}
\enq
and $R_M \in  L(X,\infty)$. Also $P(x)$ and  $\Delta_M(x)$ 
are both $o(1)$ as $x \rightarrow \infty $. We 
can now state and prove a lemma which gives the asymptotic form 
of the solutions of (1.3), incorporating an explicit estimate 
of the error term derived from a knowledge of $R_M$. In the lemma 
we write $P = (p_{ij})$ and  $\Delta_M =  {\rm{dg}} (\delta_{1M}, ..., \delta_{nM})$.
\begin{lemma}                             
  Let $M (\geq 2)$ be an integer and let $D_M$, $R_M$ and $P$ be as 
in (\ref{eq:3.8})  - (\ref{eq:3.10}). In some interval $[X,\infty)$, let
\beq
       {\rm{ Re}}\{ (\omega_ j - \omega_ k + \delta_{jM} - ë
        \delta_{kM})Q^{1/n} \} \;\;\;  (1\leq j, k \leq n)   
        \label{eq:3.11}
\enq
have constant sign (either $\geq 0$ or $\leq 0$), with
\beq
  \mid    \int_X^\infty    {\rm{ Re}}\{(\omega_ j - \omega_ k + \delta_{jM} - 
  \delta_{kM})Q^{1/n} \} dt \mid  = \infty  \label{eq:3.12}
\enq
for $j  \neq k$. Also, let
 \beq              
 \int_X^{\infty} \mid      Q^{1/n}(t)R_M(t)  dt\mid < 
 \infty.               \label{eq:3.13}
 \enq
Then, for $1 \leq k  \leq n$ and $1 \leq r \leq n$, (1.3) has solutions $y_k$ such 
that
\beq
y_k^{(r-1)} = Q^{- (n+1-2r)/2n}\{ \omega_  k^{r-1} + \sum_{j=1}^n 
\omega_j^{r-1} \bigl( \eta_j + p_{jk} +\sum_{l=1}^{n} p_{jl} \eta_l   \bigr )
\} 
  \exp \bigl (\int_X^xQ^{1/n}(t)  \{\omega_k + \delta_{kM}(t) \}dt  \bigr ) 
  \label{eq:3.14}
  \enq
where, in $[X,\infty)$, the $\eta_j$   $ \;\;(1 \leq j \leq n)$  satisfy the estimate
\beq
    \mid \eta_j(x) \mid \leq \int_X^{\infty}  
   \mid Q^{1/n}(t)\mid \parallel R_M(t)\parallel  dt/\bigl (1 - n 
    \int_{X}^\infty \mid Q^{1/n}(t)\mid \parallel  R_M(t)\parallel  dt \bigr )
    \label{eq:3.15}
    \enq
with  $\parallel R_M\parallel  = \max\  \mid r_{ijM}\mid \;\; (1 \leq j,k \leq n)$ in terms of the entries in 
$R_M$.
\end{lemma}
   {\bf{Proof}}.  By  (\ref{eq:3.11})   and  (\ref{eq:3.13})  , we can apply the Levinson 
Theorem  \cite[Theorem 1.3.1 ]{MSPE89}  to the system  (\ref{eq:3.8}), and then 
transform back to $Y$ via  (\ref{eq:3.9})   and  (\ref{eq:3.5})   to obtain vectors $Y$ 
with the form
\begin{eqnarray}
  Y &=& {\rm{dg}} (1, Q^{1/n}, ..., Q^{(n-1)/n})\Omega(I + P)(e_k + u) \nonumber \\
&\times& \exp \bigl (\int_X^x[Q^{1/n}\{ \omega_ k + \delta_{kM}\}  - (n - 1)(2n)^{-1}(Q'/Q)](t) dt \bigr ), \label{eq:3.16}
\end{eqnarray}
where $u = o(1)$ and we have used  (\ref{eq:3.10}). Let $u$ have components 
$\eta_ j\;\; (1 \leq j \leq n)$.  Then  (\ref{eq:3.14})   follows from  (\ref{eq:3.16})   on taking the $r-$
th component on each side.
  \par
 It remains to establish  (\ref{eq:3.15})   and, to do this, we proceed 
as in \cite[    (1.4.23) and (1.4.13)]{MSPE89}. It follows from  (\ref{eq:3.11})   and 
 (\ref{eq:3.12})   that, in the notation of \cite[(1.4.13) ]{MSPE89},   the non-zero
entries of $\Phi_1(x)\Phi^{-1}(t)$ 
  and 
$\Phi_2(x)\Phi^{-1}(t)$ are all $\leq 1 $ in modulus. Hence \cite[ (1.4.13)]{MSPE89}   (with $a 
= X$ here) gives
\begin{displaymath}
        \mid \eta_j(x) \mid \leq \int_X^\infty  \Biggl  (\mid Q^{1/n}r_{jkM}|  +  \sum_{l=1}^{n} \mid Q^{1/n}r_{j lM } \eta_l \mid  \Biggr ) dt.
\end{displaymath}
Thus, if $H = {\rm{max}} \mid \eta_j  \mid \; \; (X \leq x < \infty, 1 \leq j \leq n)$, we obtain
\begin{displaymath}
              H \leq \int_X^\infty \Biggl ( \mid   Q^{1/n} \mid  \parallel  R_M   \parallel + Hn \mid Q^{1/n} \mid\parallel R_M \parallel  \Biggr )dt.
\end{displaymath}
Hence     \begin{displaymath}
              H \leq \int_X^\infty \mid   Q^{1/n} \mid  \parallel  R_M   \parallel dt /\Biggl  (   1-n \int_X^\infty\mid   Q^{1/n} \mid \parallel  R_M   \parallel dt \Biggr ),
\end{displaymath}
and  (\ref{eq:3.15})   follows.

\section{ The sequence of transformations}
   In this section we define the transformations
\beq
              Z_m = (I + P_m)Z_{m+1}  \;\;\;\; (1 \leq m \leq M - 1)         
              \label{eq:4.1}
              \enq
which lead from  (\ref{eq:3.6})   to  (\ref{eq:3.8})   in the manner introduced in 
section  2(c). A typical system in the process is
\beq                     
Z'_m = Q^{1/n}(D_m + R_m)Z_m   \label{eq:4.2}             
\enq
with  (\ref{eq:3.6})   being the case $m = 1$. 
In (\ref{eq:3.6}), we emphasise the role of the diagonal terms in $C$ by taking these terms over to $D$ as in (\ref{eq:3.10}).
Thus  we  define
\beq
D_1=D+\frac{1}{2} (n-1)pI, \;\; R_1 = Q^\prime Q^{-1-1/n}(C- {\rm{dg}}C)
\label{eq:4.2a}
\enq
with $D$ as in (\ref{eq:3.7}) and $p=(Q^{-1/n})^\prime =
 -n^{-1}Q^\prime Q^{-1-1/n}$.
Then, again as in (\ref{eq:3.10}), we write 
\beq
D_m = D_1   + \Delta_ m \;\;  (m\geq 2).  \label{eq:4.3}           
\enq
\par
Now, as indicated 
by (\ref{eq:2.11}), $R_m$ will contain terms of different orders of 
magnitude as $x   \rightarrow \infty$, and we wish to identify these terms 
according to their size. Hence we write
\beq
R_m =V_m+E_m= V_{1m} + V_{2m} + ... + V_{\mu m} + E_m,        \label{eq:4.4}
\enq
where            
\beq
V_{km} = o(V_{jm})  \;\;\; (x  \rightarrow \infty , k > j)            
\label{eq:4.5}
\enq
and                 
\beq
E_m = o(V_{\mu m}) \;\;\;  (x   \rightarrow \infty).       \label{eq:4.6}          \enq
Here $E_m$ represents terms which are already of the accuracy that 
we require in (\ref{eq:3.15}), while the $V_{jm}$ represent terms which are 
not of that accuracy and which have to be replaced by smaller-order terms as we go through the transformation process in the 
manner discussed in  section 2. Also, as in the case of (\ref{eq:4.2a}), we arrange that
\beq
{\rm{dg}} V_{1m}=0.  \label{eq:4.6a}
\enq
\par
To discuss a typical step in the process leading to (\ref{eq:3.8}),
we substitute (\ref{eq:4.1}) into (\ref{eq:4.2}) to obtain   
 \begin{eqnarray}
    (I + P_m)Z^\prime_{m+1} & +& P^\prime_mZ_{m+1}\nonumber \\
           &  =& Q^{1/n}(D_m + R_m)(I + P_m)Z_{m+1} \nonumber \\
            & =& Q^{1/n}\{(I+P_m)D_m + D_mP_m -P_mD_m+ V_{1m} + V_{1m}P_m
\nonumber \\
& &+(R_m-V_{1m})(I+P_m\}Z_{m+1}.  \nonumber \\
        \label{eq:4.7}
       \end{eqnarray}
Now we define $P_m$ by
\beq
                   P_mD - DP_m = V_{1m}                 
                   \label{eq:4.8}
                   \enq
with ${\rm{dg}} P_m=0$, this definition being consistent in the diagonal entries because of (\ref{eq:4.6a}). Thus the entries $p_{ijm}$ in $P_m$ are defined by
\beq
                    p_{ijm} = v_{ij1m}/(\omega_ j - \omega_ i),  \label{eq:4.9}                       \enq
again as indicated by (\ref{eq:2.9}), (\ref{eq:2.10}) and  section 2(c). It follows from 
(\ref{eq:4.3}) and (\ref{eq:4.8}) that, in  (\ref{eq:4.7}), we have
\beq
D_mP_m - P_mD_m + V_{1m} =     
\Delta_m P_m - P_m \Delta_m =T_m, \label{eq:4.9a}
\enq
say. Hence, so far, (\ref{eq:4.7}) gives
 
\begin{eqnarray}
    Z^\prime_{m+1}& =& Q^{1/n} \{ D_{m} + (I + P_m)^{-1}
    ( -Q^{-1/n}P^\prime_m  + T_m +V_{1m}P_m\nonumber \\
    & &+ (R_m-V_{1m})  (I + P_m)) \} Z_{m+1}.    \label{eq:4.10}                                     
\end{eqnarray}
We wish to show that this can be expressed as
\beq
              Z^\prime_{m+1} = Q^{1/n}(D_{m+1} + R_{m+1})Z_{m+1},           
              \label{eq:4.7a}
 \enq
where  $R_{m+1}$ has a similar form to (\ref{eq:4.4}):
\beq 
R_{m+1} =V_{m+1}+E_{m+1} =
 V_{1,m+1} + ... + V_{\mu, m+1} + E_{m+1},  \label{eq:4.13}
\enq 
but with a different $\mu$, and where the $V_{j,m+1}$ can be obtained  
constructively from the $V_{jm}$.  Also, as in (\ref{eq:4.6a}),
 we shall arrange that
\beq
{\rm{dg}} V_{1,m+1}=0.
\label{eq:4.13a}
\enq
\par 
We  establish (\ref{eq:4.13}) by   writing 
\beq 
(I+P_m)^{-1}=I-P_m+P^2_m- ...+(-1)^\nu P^\nu_m +(-1)^{\nu+1}(I+P_m)^{-1} 
P^{\nu +1}_m. \label{eq:4.14}
\enq 
Here $\nu$ is chosen so that the product 
\begin{displaymath} 
 (I+P_m)^{-1} P^{\nu +1}_m V_{jm}, 
\end{displaymath}   
 which occurs in (\ref{eq:4.10}), has a sufficiently small order of magnitude to be  
 included with $E_m$ and form part of $E_{m+1}$. Thus   $\nu$ 
  will differ for different $j$, and similarly for the other terms 
  in (\ref{eq:4.10}).    
  We now  group  together terms of the same order of magnitude and  
  denote the dominant term by $S_{m+1}$. 
  We then    obtain (\ref{eq:4.13}) ( with $S_{m+1}$  in place of $V_{1,m+1}$),
 where $E_{m+1}$   has the same order of magnitude 
   as $E_m$ and, by (\ref{eq:4.9}), $S_{m+1}$ and the $V_{j,m+1}$ are known explicitly in terms of the  
   $V_{jm}$.
We complete the derivation of (\ref{eq:4.7a}) and ( \ref{eq:4.13}),
in which (\ref{eq:4.13a}) holds, by defining
\beq
D_{m+1} = D_m + {\rm{dg}}S_{m+1}  \label{eq:4.15a}
\enq
and 
\begin{displaymath}
V_{1,m+1}=S_{m+1}-{\rm{dg}}S_{m+1}.
\end{displaymath}
 In the next two sections, we deal with the question of  
   developing an algorithm, to be implemented in the symbolic algebra system
{ \it{Mathematica }},  for determining the $V_{jm}$ $(m=1,2,...)$ and 
   estimating the $E_m$.
   \section{The basis of the  algorithm} 
   The details of the procedure based on (\ref{eq:4.4})-(\ref{eq:4.15a}) depend on the nature of $Q$.  
   In this section we assume that 
   \beq 
   Q(x) \sim ({\rm{const.}}) x^\alpha \;\;\;\;\;\;\;(x\rightarrow \infty), \label{eq:5.1}
   \enq 
   where $\alpha >0$, and that differentiation can be carried out in the sense that 
   \beq 
   Q^{(r)}(x) =({\rm{const.}}) x^{\alpha -r}\{ 1 + o(1) \}.  \label{eq:5.2}
   \enq 
   At the end of the paper, we indicate the modifications which are made 
   to accommodate other types of $Q$. 
   \par 
   We start the detailed examination of (\ref{eq:4.4})-(\ref{eq:4.6}) 
   with the cases $m=1 $ and $2$. 
   When $m=1$, (\ref{eq:4.2}) is (\ref{eq:3.6}) and,    as in 
   (\ref{eq:4.2a}),  we have
   \beq 
   D_1=D+\frac{1}{2}(n-1)pI, \; \; \; \; R_1=Q'Q^{-1- 1/n} (C -{\rm{dg}}C)= O (x^{-a}), \label{eq:5.3}
   \enq 
   by (\ref{eq:5.1}) and (\ref{eq:5.2}), where 
   \beq 
   a=1+ \alpha/n.  \label{eq:5.4}
   \enq 
   So far, we have only $\mu=1$ and $E_1=0$ in (\ref{eq:4.4}).  Then (\ref{eq:4.10}) is simply 
   \beq 
   Z'_2=Q^{1/n}  \{ D_1 +(I+P_1)^{-1}S_2 \} Z_2, \label{eq:5.5}
   \enq 
   where 
   \beq 
   S_2 =-Q^{-1/n}P'_1+ R_1P_1. \label{eq:5.6}
   \enq 
   By (\ref{eq:4.9})
 and (\ref{eq:5.1})-(\ref{eq:5.3}), we have 
   \beq 
   P_1=O(x^{- a}),\;\;\; S_2 = O(x^{-2a}). \label{eq:5.7}
   \enq
\par 
   For the inverse matrix in (\ref{eq:5.5}), we use (\ref{eq:4.14}) and, at this point, we specify the  
   accuracy to be represented by   the $E_m$ in (\ref{eq:4.4}). We require 
   \beq 
   E_m=O(x^{-(N+2)  a}) \label{eq:5.8}
   \enq 
   for all $m$ and a fixed integer $N \geq 0$.  
   Thus the process leading to (\ref{eq:3.8}) ends when (\ref{eq:4.4}) reduces to 
   \beq 
   R_M=E_M.  \label{eq:5.9}
   \enq 
   Then by (\ref{eq:4.14}), we have 
   \begin{displaymath} 
     (I+P_1)^{-1}S_2 = \Biggr ( \sum_{j=1}^N(-1)^{j-1}P_1^{j-1}\Biggl ) S_2 
   +(-1)^N(I+P_1)^{-1}P_1^NS_2. 
   \end{displaymath} 
   Thus (\ref{eq:5.5}) is the case $m=2$ of (\ref{eq:4.2}), where (\ref{eq:4.4}) holds with $\mu=N, \; V_{j2}=(-1)^{j-1}P_1^{j-1}S_2$  $(2 \leq j \leq N)$ and 
\begin{displaymath}
E_2=(-1)^N(I+P_1)^{-1}P^N_1S_2.
\end{displaymath} 
Also, as in (\ref{eq:4.15a}),
\begin{displaymath}
D_2=D_1+{\rm{dg}}S_2,\;\; V_{12}=S_2-{\rm{dg}}S_2.
\end{displaymath}
We note that, by (\ref{eq:5.7}), 
   \beq 
   V_{j2}=O(x^{-(j+1)a}), \; \; \; E_2=O(x^{-(N+2)a})  \label{eq:5.10}
   \enq 
and, by(\ref{eq:4.3}),
\begin{displaymath}
\Delta_2={\rm{dg}}S_2=O(x^{-2a}).
\end{displaymath}
Then, since $\Delta_m$ $(m\geq3)$ is obtained from (\ref{eq:4.15a}) by adding 
successively smaller-order terms to $\Delta_2$, we have
   \beq 
   \Delta_m=O(x^{-2a})  \label{eq:5.11}
   \enq 
    for all $m$. 
 \par
We move on to general $m$ in (\ref{eq:4.4}). 
  An easy induction argument starting from (\ref{eq:5.3})
and (\ref{eq:5.10}), and based on (\ref{eq:4.9})  and (\ref{eq:4.14}),
shows that
\beq
V_{jm}=O( x^{-(m+j-1)a}). \label{eq:5.12}
\enq
Then, by (\ref{eq:5.1}) and (\ref{eq:5.11}), it follows from (\ref{eq:4.9}) and
   (\ref{eq:4.9a})  that
 \beq
P_m=O(x^{-ma}), \; \;Q^{-1/n}P'_m = O(x^{-(m+1)a}),\;\; T_m=O(x^{-(m+2)a }).\label{eq:5.13}
 \enq
It now follows from (\ref{eq:4.4}), (\ref{eq:5.12}) and (\ref{eq:5.13}) that, in ({\ref{eq:4.7a}) and (\ref{eq:4.13}),  
\beq
V_{1,m+1} = -Q^{-1/n}P'_m  +(V_{2m} -{\rm{dg}}V_{2m})  
  \label{eq:5.14}
\enq 
and then (\ref{eq:5.14}) is used in (\ref{eq:4.9}) to define $P_{m+1}$.
A further consequence of (\ref{eq:4.4}) and (\ref{eq:5.12}) is that (\ref{eq:5.9}) is achieved when
\beq
M=N+2. \label{eq:5.15}
\enq
\par
It is not difficult to check that the order relations (\ref{eq:5.12}) and 
(\ref{eq:5.13}) are in fact exact, and we can describe (\ref{eq:5.12}) by saying that $V_{jm}$ has {\it{exact order}} $m+j-1$, and similarly for (\ref{eq:5.13}).
The grouping together of terms of the same order is the basis of the algorithm
which we go on to describe in more detail now.

\section{A symbolic algorithm for the computation of  solutions  }
We can now use the ideas in section 5 to   discuss the development of a
computational  algorithm 
for computing the   solutions of   (\ref{eq:1.3})
over an interval $[X,\infty)$, subject to (\ref{eq:5.1})
(\ref{eq:5.2}).
The algorithm is implemented in the symbolic algebra system 
{\it{Mathematica}} and we use the notation in sections 3-5 to denote the symbolic objects that we need.
 The   algorithm    computes   estimates for 
$n$ linearly independent  solutions of (\ref{eq:1.3}), and for the  
derivatives, and it is structured in three distinct stages. 
In the first stage, we shall not make any assumptions about the 
order of the differential equation    while, for reasons of clarity, 
in the second and third stages our discussion is focussed on the case $n=4$
and more particularly on the equation  
\beq
y^{(4 )}-x^\alpha y=\lambda y \;\;(X\leq x < \infty ), 
\label{eq:6.1}
\enq
where $X>0$ and $\alpha >0$, this equation being covered by (\ref{eq:5.1})
and
(\ref{eq:5.2}).
\par
We recall that  the discussion in sections 4-5 of 
this paper shows how the   system
\beq
Z_m^\prime=Q^{1/n}(D_m +V_m+E_m)Z_m   \label{eq:6.0}
\enq
is transformed into
\beq
Z_{m+1}^\prime=Q^{1/n}(D_{m+1} +V_{m+1}+E_{m+1})Z_{m+1}     \label{eq:6.0a}
\enq
by the mapping
\beq
{Z_m}=(I+P_m)Z_{m+1}.               \label{eq:6.0b}
\enq
Here   $P_m$, $D_m$ and  $V_m$ are defined in (\ref{eq:4.9}), (\ref{eq:4.15a}),(\ref{eq:4.4}) and (\ref{eq:4.6a}), where $V_{jm}$
has exact order $m+j-1$ as described at the end of section 5.
We now restate some of the features of this discussion in a suitable form for implementation in our algorithm.
\par
First, the number $M$ which appears in (\ref{eq:3.8}) and (\ref{eq:5.9}) is determined by the accuracy of our working which, by (\ref{eq:5.8}) and (\ref{eq:5.15}), can now be stated as 
\beq
E_m=O(x^{-Ma})\;\; (2 \leq m \leq M).               \label{eq:6.0c}
\enq
 Then, by (\ref{eq:5.12}), we have $\mu=M-m$ in (\ref{eq:4.4}).
Next, examining the use of (\ref{eq:4.14}) in (\ref{eq:4.10}), we denote by $U$ any one of the terms in (\ref{eq:4.10}) on which the inverse acts, excepting $E_m$. Thus
\beq
U \in \{ -Q^{-1/n}P^\prime_m,T_m,V_{1m}P_m,V_{jm},V_{jm}P_m \;\;
(2 \leq j \leq M-m )\}.               \label{eq:6.0d}
\enq
For each $U$, the integer $\nu$ in (\ref{eq:4.14}) is chosen so that 
$P^{\nu+1}_mU$
has sufficiently small order of magnitude to satisfy (\ref{eq:6.0c}).
Since the order of $U$ is known from (\ref{eq:5.12}) and (\ref{eq:5.13}),
the computation of $\nu$ and the grouping of terms $(-1)^rP^r_mU$ $(0\leq r \leq \nu )$
of the same order can be performed for each $U$ and $r$.
\par
By (\ref{eq:4.4}), the error term $E_m$ appears in (\ref{eq:4.10}) in the form
\begin{displaymath}
A_mE_m(I+P_m),            
\end{displaymath}
where the symbol $A_m$ stands for $(I+P_m)^{-1}$.
This is taken as an initial definition of $E_{m+1}$ which is then updated 
as (\ref{eq:4.14}) is used for each $U$, to yield the final $E_{m+1}$ in 
( \ref{eq:4.13}).
\par
This discussion leads to the following algorithm, based on (\ref{eq:4.9}) and (\ref{eq:4.14}),
 for computing the $V_{jm}$ and $S_m$--- and hence the $D_m$ and $V_m$--- in sections 4 and 5.
  \begin{Algorithm} \label{alg:6.1}
\begin{list}%
{( \arabic{rem1} )}{\usecounter{rem1}
\setlength{\rightmargin}{\leftmargin}}
 \item
First input $M$ to fix both the number of   iterations  required for 
(\ref{eq:3.8}) and the order (\ref{eq:6.0c})  of    $E_M$. 
 \item
Start with  $D_1$ and $V_1$,  and put $E_1 =0$.
 \item
For  $m=1$ to $M-1$,
 $E_{m+1}=A_mE_{m}(I+P_m)$.
 \item
For each $U$ in  (\ref{eq:6.0d}), determine $\nu$, the number needed in (\ref{eq:4.14}).
 \item
  Update error term, i.e. $E_{m+1}=E_{m+1}+ (-1)^\nu A_m P^{\nu+1}_m U$.
  \item
For $r=0$ to $\nu$,
 determine the exact order $\rho =mr+($order of $U)$ of $P_m^r U$. 
 \item
  Update $V_{\rho-m,m+1}= V_{\rho-m,m+1} + (-1)^rP^r_m U\;\; (\rho-m \geq 2)$.
\end{list}
Then, as in (\ref{eq:5.14}), 
\begin{displaymath}
S_{m+1}=-Q^{-1/n}P^\prime_m + V_{2m} 
\end{displaymath}
and
\begin{displaymath}
V_{1,m+1}=S_{m+1} - {\rm{dg}}S_{m+1}.   
\end{displaymath}
\end{Algorithm}
\par
At any stage in the algorithm, $S_{m+1}$ depends on the terms $D_j,P_j,P_j^\prime$ and $ V_j$  $(1\leq j \leq m)$. However, simplifications can be made by retaining earlier terms $S_j$ where possible. This reduces the number of terms in the computation of $S_{m+1}$ with consequent economy in the algorithm.
We illustrate this point by giving the results for $S_3$ and $S_4$ which are first obtained in terms of $V_1, V_2, P_1$ etc, and then  the previously computed
definition of $S_2$ is used dynamically  to simplify the expressions for 
$S_3$ and $S_4$. Thus we have 
\begin{eqnarray}
S_2 &=&  - {{P^\prime_1}\over {{{Q}^{{1\over n}}}}}
+V_1 P_1\;\;\;{\rm{(as \; in \;(\ref{eq:5.6}))   } } \nonumber \\
S_3 &=& - {{P^\prime_2}\over {{{Q}^{{1\over n}}}}} + T_1  -P_1 S_2\nonumber \\
S_4 &=& - 
  {{P^\prime_3}\over {{{Q}^{{1\over n}}}}} + T_2 -P_1 T_1 + V_2 P_2 + P_1 ^2 S_2. \nonumber 
\end{eqnarray}
  We also give the error term $E_5$ as an example:
\begin{eqnarray}
E_5 &=& A_4 A_3 A_2 A_1 \left\{ {P_1}^2 T_1 - {P_1}^3 S_2 \right\} (I+P_2) (I+P_3) (I+P_4) \nonumber \\
 &+& A_4 A_3 A_2 \left\{ \left( -P_1 T_1 + {P_1}^2 S_2 - P_1 S_2 + T_1 \right) P_2 \right. \nonumber \\
& & \;\;\;\;\;\; - P_2 \left. \left( S_2 - P_1 T_1 + {P_1}^2 S_2  + T_2 + V_2 P_2 \right) \right\}  (I+P_3) (I+P_4) \nonumber \\
&+& A_4 A_3 \left\{ T_3 - P_3 S_2 + V_3 P_3 + \left( -P_1 T_1 + {P_1}^2 S_2 + V_2 P_2 + T_2 \right) P_3 \right\} (I+P_4) \nonumber \\
&+& A_4 \left( -\frac{{P_4}^{'}}{Q^{\frac{1}{n}}} + T_4 + V_4 P_4 \right). \nonumber
\end{eqnarray}
\par 
The expression    
  \begin{eqnarray} 
S_6 &=&   -Q^{-{1\over n}}P^\prime_5   + T_4+
\left( -P_1 T_1 + P_1 ^2 S_2 \right)  P_2  \nonumber \\
&-&  P_2 \left( -P_1 T_1 + V_2 P_2
+   P_1 ^2 S_2 + T_2 \right)  -  P_1 ^3 T_1  +  P_1 ^4 S_2 +V_3 P_3
\nonumber  
 \end{eqnarray}
  which we   use later
   in our examples, has an associated    
   error term $E_6$    with over 250 components. 
  The  numbers of components in 
$S_8$ and $E_8$ are about $60$ and $700$   respectively. This complexity 
indicates that the computation is intractable by hand, and that a symbolic 
algebra system is the only way to  generate the recurrences.
The algorithm is however   quite cheap to compute, and we give some sample 
times in Table \ref{tab:6.1}.
 \begin{table}[h]
\begin{center}
\begin{tabular}{||c|c||} \hline
 $m$ & duration in seconds \\ \hline
$4$ & $1.6167$ \\ \hline
$6$ & $7.5000$ \\ \hline
$7$ & $22.4167$ \\ \hline
$8$ & $61.9833$\\ \hline
\end{tabular}
\end{center}
\caption{Time in seconds for computing expressions $S_m$  and $E_m$.}
\label{tab:6.1}
\end{table}
\par
At this stage, the algorithm involves    nothing more than a series of 
recurrence formulae
in terms of objects which satisfy non-commutative multiplication and
certain order relations.   When we wish to go further and 
apply the algorithm to 
a specific equation (\ref{eq:1.3}), we need to express the recurrence  
formulae in terms of   $n \times n$ matrices  and, in the   
 next stage of the algorithm, we focus on $n=4$. 
In this case (\ref{eq:3.7}) becomes
\begin{displaymath}
 D =    
{\rm{dg}} (1 , i  , -1 ,-i  )  
\end{displaymath}
and, in (\ref{eq:4.2a}),
\beq
p=-\frac{1}{4}Q^\prime Q^{-5/4}. \label{eq:6.7m}
\enq
Also, it is easily verified from (\ref{eq:3.3}) and (\ref{eq:4.2a}) 
that $V_1\;(=R_1)$ has the form to be stated now.
  \begin{Algorithm} \label{alg:6.2}
Starting with   initial matrices
\begin{displaymath}
 D_1 =    
D+\frac{3}{2}pI, 
\end{displaymath}
\begin{displaymath}
 V_1 = -\frac{1}{2}p\left( 
\begin{array}{cccc}
0 & 1+i&1&1-i \\ 1-i&0&1+i&1\\ 1&1-i&0&1+i\\1+i&1&1-i&0
\end{array}  \right ),
\end{displaymath}
and $E_1=0$,
the expressions   $S_2, ..., S_{M-1}$  generated in Algorithm  \ref{alg:6.1} are evaluated
in order.
These are then used in turn to generate the $ 4 \times 4 $ matrices
  $D_2, ..., D_{M-1}$ ,$V_2, ..., V_{M-1}$,
$P_2, ..., P_{M-1}$,
$P^\prime_2, ..., P^\prime_{M-1}$.
\end{Algorithm}
\par
We note that the process stops at $M-1$ because (\ref{eq:5.9}) implies that
$D_M=D_{M-1}$ and   $V_M=0$. In Algorithm \ref{alg:6.2}, we assume only that
$p$ and $Q$ are 
   sufficiently smooth functions related by (\ref{eq:6.7m}) and 
 \beq
 Q''= -4Q^{\frac{5}{4}}p'+20p^2Q^{\frac{3}{2}}. \label{eq:6.2}
\enq
This equation is used to express  higher order derivatives of $Q$ in terms of 
those of $p$. 
 \par
There are several severe computational problems  contained in this seemingly 
simple process of matrix multiplication   and substitution.
 Although all the operations needed are to be found in the Mathematica system,
the use of the {\it{Substitute}} command in Mathematica requires an inordinate
amount of computational effort, and therefore it is necessary to consider what 
substitutions and simplifications are best used at each stage of the 
algorithm.
The relationships in (\ref{eq:6.7m}) and (\ref{eq:6.2}) have been found most useful in simplifying  
the expressions produced by this part of the algorithm.
 The    main problem however, manifested by long computing times, lies in the exponential growth of the number of   terms in  the elements of the matrices $S_m$.
For example, each of the elements of $S_6$ consists of a sum of  over one thousand terms, each of which is a product of several distinct   objects.
A consequence of this complexity is that, as $m$ increases, the amount of computational effort required  to perform the calculations increases exponentially. 
The time in cpu seconds needed to compute      $S_4$ and  $S_6$ on a SPARC 10 workstation,   is given in Table \ref{tab:6.2}. The - indicates that is has not been possible to complete the computation within a reasonable time. This is  because  of the limitations in  the size of memory and cpu speed that we have available.
 \begin{table}[h]
\begin{center}
\begin{tabular}{||c|c||} \hline
$ m $& duration in seconds \\ \hline
$4 $& $74.6167$ \\ \hline
$6 $&$ 904.24$ \\ \hline
$7$ &$ - $\\ \hline
\end{tabular}
\end{center}
\caption{Time in seconds for computing matrices $S_m$.}
\label{tab:6.2}
\end{table}
\par
The final stage of the  algorithm deals with the remaining matter of obtaining an
explicit upper bound for the norm $\parallel E_m \parallel $ of $E_m$
which adds precision to the order estimate (\ref{eq:6.0c}) and which can 
be used in (\ref{eq:3.15}) when $m=M$.
We use the $\sup$ norm, for which the triangle and Cauchy inequalities are
\begin{displaymath} 
\parallel A+B \parallel \leq \parallel A \parallel + \parallel B \parallel,
\;\;\; \parallel AB\parallel \leq 4 \parallel A \parallel \parallel B \parallel
\end{displaymath}
in our case of $ 4 \times 4 $ matrices. 
For the inverse matrix $A_m = ( I +P_m)^{-1} $, we have from the 
geometric series
\begin{displaymath}
\parallel A_m \parallel \leq 1+ \sum_{r=1}^\infty \parallel P^r_m \parallel \leq
1 + \sum_{r=1}^\infty 4^{r-1} \parallel P_m \parallel^r  
\end{displaymath}
and hence
\beq
\parallel A_m \parallel \leq 1 + \parallel P_m \parallel/(1-4 \parallel P_m \parallel)
\label{eq:6.n1}
\enq
provided that   $\parallel P_m \parallel <1/4$. 
\par
There is one other point to mention before discussing this final stage of the  algorithm, 
and this is prompted by  the fact that the reciprocal of $Q$ occurs throughout
 the process in Algorithms \ref{alg:6.1}  and   \ref{alg:6.2}. 
 In order to estimate this reciprocal, we require a lower bound for $Q$ in the
  relevant $x-$interval $[X, \infty)$. Thus, in addition to (\ref{eq:5.1})
  and (\ref{eq:5.2}), we require that
 \beq
 \mid Q(x) \mid \geq k^{-1} x ^ \alpha \;\;\; ( X \leq x < \infty )
 \label{eq:6.n2} 
 \enq
 for some constant $k > 1$. In (\ref{eq:6.1}) itself, 
 we have $Q(x)=\lambda + x^ \alpha$, and we shall consider 
 $\lambda$ in the half plane
 \beq
 {\rm{Re}} \lambda \geq -1.   \label{eq:6.n3} 
 \enq
 Thus $ \mid Q(x) \mid \geq \mid    {\rm{Re}} \lambda   +x^\alpha \mid
 \geq ( 1 -X^{-\alpha} ) x ^\alpha$
 provided that $X>1$. Thus (\ref{eq:6.n2})  holds with
 \begin{displaymath}
 k=(1-X^{-\alpha})^{-1}.
 \end{displaymath}
\begin{Algorithm} \label{alg:6.3}
Compute the $sup.$ norm of each matrix in $E_m$, using (\ref{eq:6.n1})
for the inverse matrices.
Next apply the triangle and Cauchy inequalities to obtain an upper bound for the $sup.$ norm of $E_m$ itself. 
 \end{Algorithm}
\par
The terms that we encounter in computing the matrix norms involve $p$ and its derivatives.  In order to obtain an upper bound for the $k^{th}$ derivative of 
$p(x)$, we first  use
\begin{displaymath}
Q(x)=\lambda + x^\alpha
\end{displaymath}
and  (\ref{eq:6.2}) to compute a symbolic expression for $p^{(k)}$. We then use  the inequality (\ref{eq:6.n2})
 to compute bounds for $\mid p^{(k)}(x)\mid $. In obtaining this estimate  
every element of every matrix  which occurs in the error term $E_M$ must be examined. The 
  triangle and Cauchy inequalities  are applied to each element and    the bounds for $\mid Q \mid $ and $\mid p^{(r)}\mid$  are substituted with  
 $ k \geq (1-X^{-\alpha})^{-1} $  and the specific values of $X$ and  $\alpha$ under consideration. 
This enables the  $sup.$ norm of each matrix to be computed  at the point $X$ and   a precise upper bound $\epsilon (X) $ for the error in the solutions 
 determined. Again the amount of computational effort needed is exponentially increasing with  the requested number of iterations.
In table \ref{tab:6.3} we show, where possible, the amount of cpu time needed to compute
 these norm estimates. 
\begin{table}[h]
\begin{center}
\begin{tabular}{||c|c|c|c||} \hline
 $m$ & \multicolumn{3}{c||}{durations in seconds} \\ \cline{2-4}
   & $S_m , R_m , P_m , A_m$ norms & $T_m$ norms & $P^\prime_m$ norms \\ \hline
$4$ & $167.833$ & $148.25 $&$ 241.33 $\\ \hline
$6 $& $1159.5 $& $1446.2$ & $2216.25$ \\ \hline
$7$ & $3297.95$ & $-$ &$ -$ \\ \hline
\end{tabular}
\end{center}
\caption{Time in seconds for computing norms where $k=\frac{11}{10}$ and $\alpha=1$.}
\label{tab:6.3}
\end{table}
\par
With the completion of the three stages of our algorithm, we have all the information needed to compute the product
\beq
 \prod_{m=1}^{M-1}(I+P_m)  \label{eq:6.3}
\enq
which appears in (\ref{eq:3.9}) and which contributes the terms $p_{ij}$ to
(\ref{eq:3.14}). Also, by (\ref{eq:5.9}),
we have the information needed for the estimate (\ref{eq:3.15})
of  $\eta_j$ which represents the error term in (\ref{eq:3.14}). Thus our algorithm computes the solutions (\ref{eq:3.14}) of (\ref{eq:1.3}) to a specified accuracy based on (\ref{eq:3.15}). We have of course chosen to orientate our discussion towards the example (\ref{eq:6.1}).
 \par
The  computation  associated with the product (\ref{eq:6.3}) can be implemented symbolically only for $M<6$ since, for larger $M$, the number of terms that need to be manipulated becomes too large   
for the computer
store that we have available,  and the computation is performed 
numerically with an accuracy of $30$ decimal digits. In table \ref{tab:6.4}
we give 
the times needed to compute (\ref{eq:3.14}) using symbolic and numerical methods.
The  shorter time needed when $\alpha=1$ is a consequence  of derivative terms,
which    are not zero for non-integer $\alpha$,   becoming   zero  when  $\alpha=1$.
\begin{table}[h]
\begin{center}
\begin{tabular}{||c|c|c|c|c|c|c||} \hline
$ M$ & \multicolumn{3}{c|}{symbolic calculation} & \multicolumn{3}{c||}{numeric calculation} \\ \cline{2-7}
   & $\alpha=\frac{1}{2}$ & $\alpha=1$ & $\alpha=\frac{4}{3}$ & $\alpha=\frac{1}{2}$ & $\alpha=1$ & $\alpha=\frac{4}{3}$ \\ \hline
$4$ & $120.65$ & $17.5167$ & $ 149.983$ & $9.55$ &$ 6.46667$ & $13.8833 $\\ \hline
$6$ & $-$ & $ 182.583$ & $ -$ & $34.1333$  & $11.2$ & $ 40.5333$ \\ \hline
$7 $&$ - $& $773.8$ &$ - $& $71.556$ & $19.35$ &$84.6667$ \\ \hline
\end{tabular}
\end{center}
\caption{Time in seconds for determining (\ref{eq:3.14}) with    $X=10$, and $\lambda = 1 + i$.}
\label{tab:6.4}
\end{table}
\section{Generalised Titchmarsh-Weyl theory}
The motivation for the asymptotic analysis and associated algorithm in sections
2-6 is provided by the spectral theory of the equation
\beq
(-1)^{\nu}y^{(2 \nu )}+q(x)y=\lambda y \;\;\; ( 0 \leq x < \infty)
\label{eq:7.1}
\enq
together with boundary conditions at $x=0$. Here $q(x)$ is real-valued and 
locally integrable in $[0, \infty)$ and $\lambda$ is a complex spectral
parameter. We take the boundary conditions to be the Dirichlet conditions
\beq
y^{(k-1)}(0)=0 \;\;\;  (1\leq k \leq \nu ) \label{eq:7.2}
\enq
but we shall comment briefly on other choices of boundary conditions 
later in this section. Our purpose here is to state 
what we need from the spectral theory of (\ref{eq:7.1}), and then 
to introduce the contribution that our methods make.
\par
The fundamental result of Titchmarsh and Weyl in \cite[Chapter 2]{ECT62} and \cite{Wey10}
concerns the case $\nu=1$ of (\ref{eq:7.1}) and ( \ref{eq:7.2}), and they 
showed that, when Im$\lambda \neq 0$, (\ref{eq:7.1}) has a non-trivial solution
$\psi(x,\lambda)$ which is $L^2(0,\infty)$. The importance of  
this solution lies  in writing $\psi$ in the form
\beq
\psi (x,\lambda ) = \theta (x,\lambda ) +m(\lambda) \phi(x,\lambda ),
\label{eq:7.3}
\enq
where $\theta$ and $\phi$ are the solutions of (\ref{eq:7.1}) ( with $\nu=1$)
which satisfy the initial conditions
\beq
\theta(0,\lambda)=1, \;\;\theta^\prime (0,\lambda)=0;\;\; 
\phi(0,\lambda)=0,\;\; \phi ^\prime (0,\lambda)=1. \label{eq:7.4}
\enq
The function $m(\lambda)$ thus defined is the basis for the spectral 
theory of  (\ref{eq:7.1}) and (\ref{eq:7.2}) as developed in \cite[Chapter 2]{ECT62} and \cite{Wey10}
for the case $\nu=1$.  Depending on the nature of $q(x)$, either $m(\lambda)$
is uniquely determined for every non-real $\lambda$ 
( the so-called limit-point case ) or $m(\lambda)$ involves an 
additional parameter   ( the limit-circle case).
We refer to the paper by Fulton \cite{CTF77} for a more recent critical discussion
of the limit-circle case. A number of explicit examples of $m(\lambda)$
are given in \cite[Chapter 4]{ECT62} for equations ( \ref{eq:7.1}) which can
be solved in terms of the special functions.
\par
The generalisation of the Titchmarsh-Weyl theory to higher-order equations, 
of which (\ref{eq:7.1}) is an instance, was made by Everitt \cite{Eve63}, \cite{Eve64}.
Here we indicate    this generalisation in the yet wider context of the
Hamiltonian system
\beq
JY^\prime (x) = \{ \lambda A(x) + B(x) \} Y(x) \;\;\;(0\leq x< \infty),
\label{eq:7.5}
\enq
the spectral theory of which was initiated by Atkinson \cite[Chapter 9]{FVA64}
and further developed by Hinton and Shaw in a series of papers which include 
\cite{HS81a}, \cite{HS81b}, \cite{HS82a}. In (\ref{eq:7.5}), $Y$ is a $2 \nu$-component vector, $A$ and $B$
are Hermitian matrices with $A \geq 0$, and
\begin{displaymath}
J=\left ( \begin{array}{cc}
0 & -I_\nu \\
I_\nu & 0 
\end{array}
\right )
\end{displaymath}
in terms of the $ \nu \times \nu$ identity matrix $I_\nu$.
We note that (\ref{eq:7.1}) is the special case of (\ref{eq:7.5})
in which the entries $a_{ij}$ and $b_{ij}$ of $A$ and $B$ are
\begin{eqnarray*}
a_{11}=1, & a_{ij}=0\;{\rm{otherwise}}, \\
b_{11}=-q, & b_{2\nu , 2 \nu}=1, 
\end{eqnarray*}
\begin{displaymath}
b_{ij}=1\;\;(j=i+\nu-1,\;2\leq i\leq \nu \;{\rm{ and }} \;i=j+\nu-1, \; 2 \leq j \leq \nu  ),
\end{displaymath}
 and $b_{ij}=0$ otherwise.
The components $y_i$ of $Y$ are then given by
\beq
y_i=y^{(i-1)}\;(1\leq i \leq \nu),\;=(-1)^iy^{(3\nu-i)}\;
(\nu+1 \leq i \leq 2 \nu ).  \label{eq:7.6}
\enq
\par
We now introduce the Dirichlet boundary condition
\beq
(I_\nu \; 0)Y(0)=0            \label{eq:7.7}
\enq
for (\ref{eq:7.5}), which corresponds to ( \ref{eq:7.2}), 
and we also introduce the $2 \nu \times \nu $ solution matrices 
$\Theta(x,\lambda)$ and $\Phi(x,\lambda)$ of (\ref{eq:7.5}) which satisfy
the initial conditions
\beq
\Theta(0,\lambda)  =
\left( \begin{array}{c} I_\nu \\0  \end{array} \right ),\;\;\;
\Phi(0,\lambda)   =\left( \begin{array}{c}0\\ I_\nu  \end{array} \right
).
\label{eq:7.8}
\enq
These conditions correspond to (\ref{eq:7.4}). In the references just cited,
it is shown that, when Im$\lambda \neq0$, (\ref{eq:7.5}) has a $2\nu \times 
\nu$ solution matrix, with rank $\nu$, such that 
\beq
\int_0^\infty \Psi^*(x,\lambda)A(x)  \Psi(x,\lambda)dx < \infty.
\label{eq:7.9}
\enq
Further, corresponding to (\ref{eq:7.3}), $\Psi$ can be written as
\beq
\Psi(x,\lambda)=\Theta  (x,\lambda)   + \Phi (x,\lambda) M(\lambda),
\label{eq:7.10}
\enq
so defining the $\nu \times \nu$ matrix $M(\lambda)$.
This matrix again forms the basis for the spectral theory of (\ref{eq:7.5})
and (\ref{eq:7.7}), as in the case of $m(\lambda)$ in (\ref{eq:7.3}).
\par
In analogy with the Titchmarsh-Weyl theory, we say that (\ref{eq:7.5})
is in the {\it{limit-point case}} if, for some $\lambda$ with 
Im$\lambda \neq0$, (\ref{eq:7.5}) has no more that $\nu$ linearly independent
 solutions such that 
 \begin{displaymath}
 \int_0^\infty Y^*(x,\lambda) A(x)Y(x,\lambda) dx < \infty.
 \end{displaymath}
 This limit-point classification is independent of $\lambda$ and,
 by (\ref{eq:7.9})
 and (\ref{eq:7.10}), it has the implication that $M(\lambda)$ is uniquely
 determined for all non-real $\lambda$. Further, if $\Psi_1(x,\lambda)$ is
 any other $2 \nu \times \nu$ solution matrix of (\ref{eq:7.5}) which has
 rank $\nu$ and also satisfies (\ref{eq:7.9}), then
 \beq
 \Psi_1(x,\lambda)=\Psi(x,\lambda)C(\lambda) \label{eq:7.11}
 \enq
 for some non-singular $\nu \times \nu $ matrix $C$. It now follows from 
 (\ref{eq:7.8}),(\ref{eq:7.10}) and (\ref{eq:7.11})  that
 \begin{displaymath}
 \Psi_1(0,\lambda) =\left( \begin{array}{c} C(\lambda) \\M(\lambda) C(\lambda) ) \end{array} \right ).
 \end{displaymath}
Hence, partitioning
 \beq
 \Psi_1(0,\lambda) =\left( \begin{array}{c} \zeta(\lambda) \\\eta(\lambda)   \end{array} \right ),  \label{eq:7.12}
 \enq
we obtain
 \begin{displaymath}
\zeta(\lambda) =C(\lambda), \;\;\;\eta(\lambda)=M(\lambda)C(\lambda).  
 \end{displaymath}
 Eliminating $C(\lambda)$, we obtain a formulae for $M(\lambda)$ in terms of the initial values of $\Psi_1$ at $x=0$:
\beq
M(\lambda)=\eta(\lambda)\zeta^{-1}(\lambda).   \label{eq:7.13}
\enq
Finally here, we note that $M(\lambda)$ is analytic for Im$\lambda \neq0$
and $M(\lambda)$ coincides with its own transpose, that is, its entries satisfy
$m_{ij}(\lambda)=m_{ji}(\lambda)$  \cite{HS81b}.
\par
For use in the next section, we mention here the form which (\ref{eq:7.13})
takes when (\ref{eq:7.5}) arises from the fourth-order example of (\ref{eq:7.1})
with $\nu=2$. In this situation, we have two linearly independent $L^2(0,\infty)$ solutions    $\psi_1 (x,\lambda)$ and $\psi_2(x,\lambda)$
of (\ref{eq:7.1}). Then, by (\ref{eq:7.6}), (\ref{eq:7.12})
and (\ref{eq:7.13}), the spectral matrix $M(\lambda)=\bigl(m_{ij}(\lambda)\bigr)$ is given by
\beq
\bigl(m_{ij}(\lambda)\bigr)=
\left(
\begin{array}{cc}
-\psi^{\prime\prime\prime}_1(0,\lambda) & 
-\psi^{\prime\prime\prime}_2(0,\lambda) \\
\psi^{ \prime\prime}_1(0,\lambda) & \psi^{ \prime\prime}_2(0,\lambda) 
\end{array}
\right)
\left(
\begin{array}{cc}
\psi_1(0,\lambda) & \psi_2(0,\lambda) \\
\psi^{ \prime}_1(0,\lambda) & \psi^{ \prime}_2(0,\lambda) 
\end{array}
\right)^{-1}. \label{eq:7.14}
\enq
\par
We have mentioned that examples of $m(\lambda)$ exist when $\nu=1$
in (\ref{eq:7.1}), arising from the special functions. In contrast, no such examples are known when $\nu \geq 2$ except for the relatively trivial Fourier case when $q(x) =0$.
In order therefore to illustrate and test the general spectral theory of ( \ref{eq:7.1}), there is a need for a computational approach to estimate $M(\lambda)$ numerically. Thus the motivation for our work
is the development of an effective computational approach and, in the next section, we give our numerical results for the example where $\nu=2$ and $q(x)=-x^{\alpha}$. Our methods are not confined to this choice of $q(x)$, but this choice does have a certain significance which we discuss in section 10.
\par
We have concentrated on the Dirichlet boundary conditions (\ref{eq:7.2})
and (\ref{eq:7.7}). If these are changed to other conditions which involve linear combinations of derivatives, then (\ref{eq:7.8}) and consequently
$M(\lambda)$ are also changed. However, in the limit-point case emphasised here, the new $M(\lambda)$ is related to the Dirichlet $M_D(\lambda)$ by a
functional equation, and therefore in principle the properties of the new $M(\lambda)$ can be obtained from those of $M_D(\lambda)$. We refer to \cite[p.69]{ECT62} and \cite[Theorem 3.4]{NGJD} for the details of the functional equation. Here, as an example, we note only that, if (\ref{eq:7.7}) is
replaced by the Neumann condition $(0 I_\nu)Y(0)=0$,
then \cite[Corollary 3.1]{NGJD}
\begin{displaymath}
M_N(\lambda)=-M_D^{-1}(\lambda).
\end{displaymath}

\section{Computation of the spectral matrix}
We can now draw together the ideas in the previous sections to give an effective procedure for the numerical computation of the spectral matrix $\bigl( m_{ij}(\lambda) \bigr)$ associated with the equation
\beq
y^{(4)} - x^\alpha y = \lambda y \;\;\; ( 0 \leq x < \infty) \label{eq:8.1}
\enq
and the Dirichlet boundary conditions
\begin{displaymath}
y(0)=0, \;\;y^\prime(0)=0. 
\end{displaymath}
The range for $\alpha$ is $ 0  < \alpha \leq 4/3$, so that ( \ref{eq:8.1})
is in the limit-point case as required in section 7 \cite[Theorem 3.11.1(b)]{MSPE89} and \cite[section 23.4]{MAN68}.
\par
The $m_{ij}(\lambda)$ are given by ( \ref{eq:7.14}) and we therefore require the values  at $x=0$ of two $L^2(0,\infty)$ solutions $\psi_1(x,\lambda)$ and
 $\psi_2(x,\lambda)$ of (\ref{eq:8.1}), together with the values of their 
 derivatives. We take Im$\lambda>0$ to be definite. 
 Then consideration of the exponential factor in (\ref{eq:3.14}), 
 with $Q=\lambda+x^\alpha$ and $n=4$, shows that the two $L^2(0,\infty)$ 
 solutions of (\ref{eq:8.1}) are associated with $\omega_2$ and  $\omega_3$. 
 Thus (\ref{eq:3.14}) provides the values of $\psi_1(x,\lambda)$ and 
 $\psi_2(x,\lambda)$ and their derivatives at $x=X$, with an error given 
 by (\ref{eq:3.15}). 
These values are produced by our algorithm as explained at the end of section 6,
and we denote
 the maximum error in this computation by $\epsilon(X)$.
The values of the solutions and their derivatives at $X$ are used as 
initial data for a numerical initial value solver 
which is then used to compute numerically the values of the solutions 
and their derivatives at $0$.  These values are then used in (\ref{eq:7.14}).
This computation back to $x=0$ has been performed with the NAG library code D02NMF in the following way.
We introduce the standard Riccati variable
$\xi(x)= \sigma(x) \tau^{-1}(x)$, where
\begin{displaymath}
\sigma = \left ( \begin{array}{cc}
\psi_1 & \psi_2 \\
\psi_1^\prime & \psi_2^\prime
\end{array}
\right ),
\;\;
\tau = \left ( \begin{array}{cc}
-\psi_1^{\prime \prime \prime } & -\psi_2^{\prime \prime \prime } \\
\psi_1^{  \prime \prime } & \psi_2^{  \prime \prime }
\end{array}
 \right ).
\end{displaymath}
Then, as in \cite{BBEMM95}, $\xi$ satisfies  a first-order non-linear differential equation. The above code is applied to $\xi$, and it computes $\xi$
back from $X$ to $0$. Then, by (\ref{eq:7.14}), $M(\lambda)= \xi^{-1}(0)$. 
\par
The value of $X$  is at our disposal, and the actual choice is governed by two factors. A smaller choice has the advantage that the integration back from $X$ to zero is performed over a smaller range with a better retention of the accuracy represented by $\epsilon(X)$. On the other hand, a larger choice has the advantage that $\epsilon (X)$ is smaller. The first  factor appears to be more significant and we have chosen $X=10$. Later we shall comment on $X$=20 and on the very large values of $X$ that arise in the alternative method of Bennewitz et al \cite{BBEMM95}. The value of $M$---the number of requested transformations in sections 3 to 6---is also at our disposal, and we have chosen $M=6$ to obtain the good accuracy which we present later in Tables \ref{tab:8.1}
-\ref{tab:8.2}
 We turn now to discuss the results of numerical experiments performed with our algorithm. In Tables \ref{tab:8.1}
-\ref{tab:8.2}
 we give the values of $m_{11}(\lambda)$, $m_{22}(\lambda)$ and $m_{12}(\lambda)=m_{21}(\lambda)$ for three values of $\alpha$ and a selection of values of $\lambda$. The values of $\alpha$ are, first, $\alpha=1$ because in this case we can obtain independent confirmation of our results from the theory of the higher-order Airy equation, which is the subject of section 9. Then $\alpha = 4/3$ is chosen because it is the maximum value of $\alpha$ for which 
(\ref{eq:8.1}) remains in the limit-point case. Finally $\alpha=1/2$ is chosen in order to give a spread of values of $\alpha$ in the range $0 < \alpha \leq 4/3$. 
\par
The values of $\lambda$, all with Im$\lambda>0$, are mainly taken to be near the origin because we wish to compare the power of our method with that of the code developed previously by Bennewitz et al. \cite{BBEMM95}. Our method is considerably more effective, and the improvement is most conspicuous for smaller values of $\mid \lambda \mid $, particularly when $\alpha =4/3$, as can be seen from Table \ref{tab:8.3} and the choice of $X$ in the other tables. We recall that a value of $X$ as small as possible is generally the most desirable because
 of the integration of (\ref{eq:8.1}) over $[0,X]$ starting from initial values at $X$.
  \subsection{Example $\alpha=1$}
We have carried out the computation using $M=6$ in the algorithm in section 6. With $X=10$ and $X=20$, we obtain
\begin{displaymath}
\epsilon(10)=2.036696 \times 10^{-6}, \;\;\; \epsilon(20)=1.038152 \times 10^{-8}.
\end{displaymath}
The numerical integration over $[0,X]$ is performed with an accuracy of $10^{-10}$. Eight values of $\lambda$ are considered, as listed in Table
\ref{tab:8.1}. In the third column, the agreement in terms of significant figures is given for the choice $X=20$. The slight loss of accuracy is in line with  the comments which we made earlier in this section concerning the choice of $X$. The last two columns give the corresponding performance of the code of Bennewitz et al \cite{BBEMM95}. The large values of $X$ and the cases of failure will be noted, as will also the cases of agreement.
\begin{table}
\small
\begin{center}
\begin{tabular}{||c|c|c|c|c||} \hline
$\lambda$ & \multicolumn{2}{c|}{Asymptotic solution}   & \multicolumn{2}{c||}{Bennewitz} \\ 
	  & \multicolumn{2}{c|}{}           & \multicolumn{2}{c||}{} \\ \hline
	  & \multicolumn{2}{c|}{numerical integration}       & accuracy & integration \\ 
	  & over $[0,10]$       & over $[0,20]$               &   & distance $X$ \\ \hline

 i       & -0.1844321489 + 1.8220405579 i & 10 sf   & 10 sf & 11755.07 \\ \hline
	& -0.6613801122 + 1.1030735970 i & 10 sf   & 10 sf &  \\ \hline
	& -1.4291142225 + 0.8401392698 i & 10 sf   & 10 sf &  \\ \hline
\hline
0.5 + i & 0.1935073882 + 2.0560526848 i & 8 sf   & Failed & 74606.02 \\ \hline
		  & -0.5871452689 + 1.2561461926 i   & 8 sf &  & \\ \hline
		  & -1.4407042265 + 0.9327853918 i   & 8 sf &  & \\ \hline
\hline
 10 + 10 i       &  2.0868480206 + 10.2971620560 i  & 10 & 10 sf & 59.15 \\ \hline
		&  -1.4532939196 + 3.5441000462 i   & 10 & 10 sf & \\ \hline
		&  -2.3046081066 + 1.5476453304 i   & 10 & 10 sf & \\ \hline
\hline
 1 + 0.001 i & 1.0726006031 + 1.6053582430 i & 7 sf   &  & Failed \\ \hline
		    & -0.2162801623 + 1.3155398369 i & 7 sf   &  & \\ \hline
		    & -1.2967267036 + 1.0783087015 i & 7 sf   &  &  \\ \hline
\hline
0.01 i & 0.3271814883  + 1.0330723524 i & 7 sf   &   & Failed \\ \hline
		    & -1.2149009705 + 0.8802024126 i & 10 sf   &  & \\ \hline
		    & -0.3127440214 + 0.9515559077 i & 7 sf   &  &  \\ \hline
\hline
\end{tabular}
\end{center}
\caption{   $\alpha=1$ }
\label{tab:8.1}
\end{table}
\subsection{Example $\alpha=4/3$}
We have again carried out the computations using $M=6$, and we obtain
\begin{displaymath}
\epsilon(10)=6.57898 \times 10^{-6}.
\end{displaymath}
Table \ref{tab:8.3} gives the values of $m_{ij}$ for eight values of $\lambda$ as before and the choice $X=10$. In the last two columns of the table, the poor performance of the previous code \cite{BBEMM95} will be noted.
\begin{table}
\begin{center}
\begin{tabular}{||c|c|c|c||} \hline
$\lambda$ & Asymptotic solution & \multicolumn{2}{c||}{Bennewitz} \\ 
	  &  & \multicolumn{2}{c||}{} \\ \hline
	  &         &  accuracy & integration \\
	  &         &     & distance $X$ \\ \hline

 i       & -0.2790325582 + 1.8323447704 i &  & Failed \\ \hline
	& -0.7153936028 + 1.1025729179 i &  &  \\ \hline
	& -1.4712435007 + 0.8293465972 i &  &  \\ \hline
\hline
0.5 + i &  0.0974698886 + 2.0557093620 i &  & Failed \\ \hline
	& -0.6384261250 + 1.2467917204 i &  &  \\ \hline
	& -1.4768083096 + 0.9157238603 i &  &  \\ \hline
\hline
 10 + 10 i &  2.0430564880 + 10.2711343765 i & 6 sf & 290.36 \\ \hline
	  & -1.4629243612 + 3.5318126678 i & 6 sf & \\ \hline
	  & -2.3070857525 + 1.5415457487 i & 6 sf & \\ \hline
\hline
 1 + 0.001 i &  0.9584534168 + 1.5867842436 i &  & Failed \\ \hline
	    & -0.2741018832 + 1.2855821848 i &  & \\ \hline
	    & -1.3294038773 + 1.0418083668 i &  &  \\ \hline
\hline
  0.01 i &  0.2010058761 + 1.0459299088 i &  & Failed \\ \hline
	    &  -1.2738007307 + 0.8492208123 i &  & \\ \hline
	    & -0.3926285803 + 0.9405522943 i &  &  \\ \hline
\hline

\end{tabular}
\end{center}
\caption{    $\alpha = \frac{4}{3}$  }
\label{tab:8.3}
\end{table}
\subsection{Example $\alpha=1/2$ } 
Again with $M=6$,we obtain
\begin{displaymath}
\epsilon(10)=1.452392 \times 10^{-6}
\end{displaymath}
and the values of the $m_{ij}$ are given in Table \ref{tab:8.2}. This time the 
previous code \cite{BBEMM95} performs better but it still does not match ours.
  \begin{table}
\begin{center}
\begin{tabular}{||c|c|c|c||} \hline
$\lambda$ & Asymptotic solution & \multicolumn{2}{c||}{Bennewitz} \\ 
	  &  & \multicolumn{2}{c||}{} \\ \hline
	  &         &  accuracy & integration \\
	  &         &     & distance $X$ \\ \hline

 i       &  0.0091688652 + 1.8079411983 i & 7 sf &  174.35 \\ \hline
	& -0.5696548223 + 1.1018460989 i & 7 sf &  \\ \hline
	& -1.3599216938 + 0.8523176312 i & 10 sf &  \\ \hline
\hline
0.5 + i &  0.3899112046 + 2.0596635342 i & 6 sf & 193.55 \\ \hline
	& -0.4997293651 + 1.2700277567 i & 7 sf &  \\ \hline
	& -1.3820117712 + 0.9563996792 i & 7 sf &  \\ \hline
\hline
 
10 + 10 i &  2.2008070946 + 10.3476095200 i & 6 sf  & 33.84 \\ \hline
	  & -1.4308696985 + 3.5664336681 i & 6 sf & \\ \hline
	  & -2.2991588116 + 1.5578293800 i & 7 sf & \\ \hline
\hline
 1 + 0.001 i &  1.2971330881 + 1.6384627819 i &  & Failed \\ \hline
	    & -0.1199851707 + 1.3626530170 i &  & \\ \hline
	    & -1.2458260059 + 1.1335541010 i &  &  \\ \hline
\hline
0.01 i &  0.5758281350 + 1.0167049170 i & - & 243897.42 \\ \hline
	    & -1.1166948080 + 0.9334350824 i & - & \\ \hline
	    &  -0.1757378578 + 0.9719497561 i & - &  \\ \hline
\hline

\end{tabular}
\end{center}
\caption{   $\alpha = \frac{1}{2}$  }
\label{tab:8.2}
\end{table}
It will be noted, in all three numerical examples that we discuss,
that  our methods produce results that are in agreement with the results in 
\cite{BBEMM95}.
  Further the algorithm reported on in       \cite{BBEMM95} fails to compute the spectral matrix for values of $\lambda$ 
close to the real line. Our algorithm has no difficulty computing the spectral matrix at these values. 

\section{The higher-order Airy equation}
Independent confirmation of our numerical results   for the case
$\alpha=1$ in Example 8.1 is provided by the theory of the equation
\beq
\frac{d^ny}{dz^n} = (-1)^n z y   \label{eq:A.1}
\enq
which, when $n=2$, becomes the well-known Airy equation. We note that 
(\ref{eq:A.1}) is itself a special case of the equation
\begin{displaymath}
\frac{d^ny}{dz^n} = (-1)^n z^m y ,
\end{displaymath}
with $m$ rational, which has been extensively investigated by Turrittin
\cite{HLT50}, Heading  \cite{JH57}, and others ( see also Paris and Wood  \cite[ pp. 188-190] {PW86}).
Here $z$ is a complex variable and, when
\beq
z=\lambda + x  \label{eq:A.2}
\enq
and $n=4$, ( \ref{eq:A.1}) becomes
\beq
y^{(4)}(x) =(\lambda +x) y(x), \label{eq:A.3}
\enq
which is the case $\alpha=1$ of (\ref{eq:8.1}).  Thus we consider the equation
\beq
\frac{d^4y}{dz^4} =   z y \label{eq:A.4}
\enq
and our primary source for the nature of the solutions is the paper by Heading \cite[ sections 1-5]{JH57}. 
\par
We note first that a straightforward power series substitution for $y$ in (\ref{eq:A.4}) gives the general solution in the form
\beq
y(z)=\sum_{r=0}^3 c_r z^rf_r(z), \label{eq:A.5}
\enq
where the $c_r$ are constants and 
\begin{eqnarray}
f_0(z) &=& 1 + \frac{1!}{5!} z^5+ \frac{1!6!}{5!10!} z^{10}
+ \frac{1!6!11!}{5!10!15!} z^{15} + ... \nonumber \\
f_1(z) &=& 1 + \frac{2!}{6!} z^5+ \frac{2!7!}{6!11!} z^{10}
+ \frac{2!7!12!}{6!11!16!} z^{15} + ... \nonumber \\
f_2(z) &=& 1 + \frac{3!}{7!} z^5+ \frac{3!8!}{7!12!} z^{10}
+ \frac{3!8!13!}{7!12!17!} z^{15} + ... \nonumber \\
f_3(z) &=& 1 + \frac{4!}{8!} z^5+ \frac{4!9!}{8!13!} z^{10}
+ \frac{4!9!14!}{8!13!18!} z^{15} + ..., \label{eq:A.6}
\end{eqnarray}
and the power series converge for all $z$. As mentioned by Heading  
\cite[p.405]{JH57}, the $f_r$ can if necessary be expressed in terms of the generalised hypergeometric function $ _0F_3$: for example
\begin{displaymath}
f_0(z) = \;_0F_3(4/5,3/5,2/5;z^5/5^4).
\end{displaymath}
\par
A particular choice of the constants $c_r$ is identified by Heading as producing solutions which have specific asymptotic behaviours for large $\mid z \mid$.
The choice \cite[(37)]{JH57} is
\begin{eqnarray}
c_0  &=&\Gamma(\frac{1}{5}) \Gamma(\frac{2}{5})  \Gamma(\frac{3}{5})   \nonumber \\
c_1  &=&\Gamma(-\frac{1}{5}) \Gamma(\frac{1}{5})  \Gamma(\frac{2}{5})5^{-4/5}    \nonumber \\
c_2   &=&\Gamma(-\frac{2}{5}) \Gamma(-\frac{1}{5})  \Gamma(\frac{1}{5})5^{-8/5} \nonumber \\
c_3   &=&\Gamma(-\frac{3}{5}) \Gamma(-\frac{2}{5}).  \Gamma(-\frac{1}{5})5^{-12/5}. \label{eq:A.7}
\end{eqnarray}
Then, as in \cite[(24),(37),(38) ]{JH57},  we denote the solution (\ref{eq:A.5}) with this choice by $J(z)$. It is shown in \cite[ (39)]{JH57}  that
\beq
J(z) \sim ( {\rm{const.}}) z^{-3/8} \exp( -\frac{4}{5}z^{5/4})
\label{eq:A.8}
\enq
for large $\mid z \mid $ with  arg $z\neq \pi$, and we can use this asymptotic formula to find the two $L^2(0,\infty)$ solutions $\psi_1(x,\lambda)$ and $\psi_2(x,\lambda)$ of (\ref{eq:A.3}) when  Im $\lambda \neq 0$.
\par
First, with $z=x+\lambda$ as in (\ref{eq:A.2}), it is clear from (\ref{eq:A.8})
 that $J(x+\lambda)$ is exponentially small as $ x \rightarrow \infty$ and hence is $L^2(0,\infty)$. Thus
\beq
\psi_1(x,\lambda)=J(x+\lambda) = \sum_{r=0}^3c_r(x+\lambda)^rf_r(x+\lambda),
\label{eq:A.9}
\enq
where the $f_r$ and $c_r$ are defined by (\ref{eq:A.6}) and (\ref{eq:A.7}).
\par
Second, as observed in \cite[(25)]{JH57},  $y(\omega z)$ is also a solution of (\ref{eq:A.4})
whenever $y(z)$ is a solution and $\omega^5=1$. We therefore consider
\beq
\psi_2(x,\lambda)=J\{(x+\lambda)e^{- 2 \pi i/5} \} = \sum_{r=0}^3c_re^{- 2 \pi ir/5}(x+\lambda)^rf_r(x+\lambda).
\label{eq:A.10}
\enq
For this solution, the exponential factor in (\ref{eq:A.8}) is 
\begin{displaymath}
\exp \{ \frac{4}{5} i ( x + \lambda )^{\frac{5}{4}} \} =
\exp \{ \frac{4}{5} i  x  ^{\frac{5}{4}}( 1 +\frac{5}{4} \lambda x ^{-1}
+O(x^{-2}))\},
\end{displaymath}
and the  modulus of this factor is
\begin{displaymath}
 \{1 +O(x^{-\frac{3}{4}})\} \exp \{ -({\rm{Im} }\lambda )  x  ^{\frac{1}{4}}\}.
\end{displaymath}
Hence, when Im$\lambda >0$,
 $\psi_2(x,\lambda)$ is also exponentially small as $ x \rightarrow \infty $ and is therefore $L^2(0,\infty)$.
\par
The initial values of $\psi_1$ and $\psi_2$ at $x=0$ are now known in 
terms of $\lambda$ from (\ref{eq:A.9})  and (\ref{eq:A.10}),
and therefore (\ref{eq:7.14}) provides an explicit formula for   $M(\lambda)$ 
in terms of   power series involving $\lambda$. In using this formula to verify our computations in Example 8.1, we have 
  approximated $\psi_1(x,\lambda)$ and $\psi_2(x,\lambda)$  by using the first $20$ terms
of the power series given in (\ref{eq:A.5}) and (\ref{eq:A.6}).
Expressions for the derivatives  have been calculated symbolically,
and the evaluation of $\psi_1(0,\lambda)$,   $\psi_2(0,\lambda)$  and the derivatives has been performed with 
$30$ digits of accuracy.
   \section{Concluding remarks}
\begin{list}%
{( \alph{rem1} )}{\usecounter{rem1}
\setlength{\rightmargin}{\leftmargin}}
 \item
{\bf{The algorithm of Bennewitz et al. \cite{BBEMM95}.}} We have already referred to some numerical results of the algorithm \cite{BBEMM95} in section 8 when $q(x)=-x^\alpha$. A more general comment on the algorithm is that
it works well in cases where $M(\lambda)$ is meromorphic, as for example when $q(x)=x^\alpha$ ($\alpha >0$). However, the computational complexity of the method causes the performance to  degrade and even fail for examples of (\ref{eq:7.1}) when arg$\lambda$ is small and $M(\lambda)$ is not meromorphic.
It is partly for this reason that we have concentrated on $q(x)=-x^\alpha$,
($0 <\alpha \leq 4/3$), the spectrum in this case filling the whole real axis
\cite[section 24.4]{MAN68} with, consequently, $M(\lambda)$
being non-meromorphic. This is the wider context within which to place our comments in section 8 concerning the effectiveness of our methods as compared to those of \cite{BBEMM95}.
\item
{\bf{HELP inequalities.}}  In the case $\nu=1$ of (\ref{eq:7.1}), the behaviour
of $m(\lambda)$ in the neighbourhood of $\lambda=0$ determines the validity of what is known as the HELP ( Hardy-Everitt-Littlewood-Polya) inequality. We refer to \cite{BEE94} and \cite{EE82} for surveys and an extensive bibliography concerning the inequality. An extension to fourth-order differential equations was given by Russell \cite{Rus78} \cite{Rus79}, but recently a more systematic development for (\ref{eq:7.1}) with general $\nu$ has been given by Dias \cite{NGJD}. Again the validity of an inequality
 of the HELP type depends on the behaviour of $M(\lambda)$ near to $\lambda=0$.
Having in this paper
 developed a computational algorithm which is effective for small $\mid \lambda \mid $, we propose to investigate further the application to the validity of higher-order HELP inequalities.
\item
{\bf{Coefficients with an oscillatory factor.}}
It is pointed out in \cite[Example 2.4.1]{MSPE89} that the method of repeated diagonalization for (\ref{eq:1.2}) works not only for coefficients of the type (\ref{eq:5.1}) but also for coefficients such as
\begin{displaymath}
Q(x)=x^\alpha p( x^\beta)
\end{displaymath}
with $\beta >0$ and $p(t)$ periodic in $t$ and nowhere zero. It seems likely that, with suitable modifications to the algorithm in section 6, the methods of this paper will cover such a coefficient in (\ref{eq:1.3}). It would appear
that a higher value of $M$ than (\ref{eq:5.15}) is needed to achieve the same accuracy (\ref{eq:5.8}) as before and that the grouping of terms in the algorithm of section 5 and 6 depends on the value  of $\beta$. These details
are  another matter for further investigation.
\item
{\bf{Hamiltonian systems in general.}}
We have concentrated in this paper on the equations (\ref{eq:1.3}) and,
for the spectral theory, (\ref{eq:7.1}).
However, our methods are in principle applicable to other differential equations and indeed Hamiltonian systems which, after an initial transformation, can be written in the form (\ref{eq:3.1}) with $\Lambda$ similar to (\ref{eq:3.2}). Such a $\Lambda$ would be
\begin{displaymath}
\Lambda = \rho (x) {\rm{dg}} ( d_1,...,d_n),
\end{displaymath}
where there is a factor $\rho (x)$ and the $d_k$ are distinct non-zero constants. More difficult however---and this is the point of this subsection---are systems (\ref{eq:3.1}) where the diagonal terms in $\Lambda$ have different orders of magnitude as $x \rightarrow \infty$. Such a situation can occur for example with the equation
\begin{displaymath}
y^{(4)}(x) +\{ P(x)y^\prime(x)\}^\prime + Q(x) y(x) =0,
\end{displaymath}
where the middle coefficient is dominant in the sense that $Q=o(P^2)$
as $x\rightarrow \infty$ \cite[section 3.5]{MSPE89}. Thus a further stage in the development of our asymptotic analysis and the associated algorithm would be to cope with differential equations, and  more generally Hamiltonian systems, where this type of $\Lambda$ occurs.
\item
{\bf{Automatic differentiation}}
The algorithm discussed in this paper has been developed using both {\it{Mathematica}} and {\it{Fortran77}}. This has resulted in  some computationally intensive symbolic  calculations as is demonstrated by the results in Tables \ref{tab:6.2} and \ref{tab:6.3}. A possible alternative
approach to performing the computation would be to use an automatic differentiation package to evaluate numerically  the  required derivatives.
Such an approach would mean that the complete algorithm could be implemented in, say, {\it{Fortran 77}} or  {\it{Fortran 90}}, thus avoiding both the long symbolic calculations reported on in section 6, and also  the need to interface the Fortran code to the {\it{Mathematica}} package. We  intend to address this issue in a future publication.
 \item
{\bf{Provably correct computations}}
There is a considerable interest in performing provably correct computations.
We note that the symbolic algorithm described in section 6 provides, not only an estimate of the $L^2[0,\infty)$ solution at some point $X>0$, but also a provably
correct bound on    the error at $X$.   This information could be used as input data to an interval-based  ordinary differential equation solver,
and the computation of the $M(\lambda)$ spectral matrix could be performed in a provably correct manner, giving precise information on the numerical errors involved. The effectiveness of this approach will need further investigation.
\end{list}

\section*{Acknowledgement}
The authors wish to thanks the EPSRC for grant GR/J61442 under which this research has been supported.

\bibliographystyle{plain}
\bibliography{biblio}  
  \end{document}